\documentstyle[amscd,amssymb,verbatim,epsf,mathrsfs,graphicx,hyperref]{amsart}

\begin{document}
\theoremstyle{plain}
\newtheorem{Thm}{Theorem}
\newtheorem{Cor}{Corollary}
\newtheorem{Ex}{Example}
\newtheorem{Con}{Conjecture}
\newtheorem{Main}{Main Theorem}
\newtheorem{Lem}{Lemma}
\newtheorem{Prop}{Proposition}

\theoremstyle{definition}
\newtheorem{Def}{Definition}
\newtheorem{Note}{Note}

\theoremstyle{remark}
\newtheorem{notation}{Notation}
\renewcommand{\thenotation}{}

\errorcontextlines=0
\renewcommand{\rm}{\normalshape}%

\title[A Converging Lagrangian Flow]%
   {A Converging Lagrangian Flow in the \\ Space of Oriented Lines}

\author{Brendan Guilfoyle}
\address{Brendan Guilfoyle\\
          Department of Computing and Mathematics  \\
          Institute of Technology, Tralee \\
          Clash \\
          Tralee  \\
          Co. Kerry \\
          Ireland.}
\email{brendan.guilfoyle@@ittralee.ie}
\author{Wilhelm Klingenberg}
\address{Wilhelm Klingenberg\\
 Department of Mathematical Sciences\\
 University of Durham\\
 Durham DH1 3LE\\
 United Kingdom.}
\email{wilhelm.klingenberg@@durham.ac.uk }

\date{18th June 2015}

\begin{abstract}
Under mean radius of curvature flow, a closed convex surface in Euclidean space is known to expand exponentially to infinity. 
In the 3-dimensional case we prove that the oriented normals to the flowing surface converge to the oriented normals of a 
round sphere whose centre is determined by the initial surface.

To prove this we show that the oriented normal lines, considered as a surface in the space of all oriented lines, 
evolve by a parabolic flow which preserves the Lagrangian condition. Moreover, this flow converges to a holomorphic Lagrangian section, 
which form the set of oriented lines through a point. 

The coordinates of this centre point are projections of the support function into the first non-zero eigenspace of the spherical Laplacian and are
given by explicit integrals of initial surface data.

\end{abstract}
\thanks{A video of a talk explaining the methods and results of this paper can be found at the following link: \url{https://youtu.be/Sx1T6legtgQ}}
\maketitle

Consider the evolution of a sphere $f:S^n\times[0,\infty)\rightarrow {\mathbb R}^{n+1}$ by {\it mean radius of curvature flow} (MRCF):
\begin{equation}\label{flow}
\frac{\partial f}{\partial t}^\perp=\sum_{j=1}^{n}r_j\;{\mbox{{\bf N}}},
\end{equation}
where ${\mbox{{\bf N}}}$ is the unit normal vector and $r_1,r_2,...,r_n$ are the radii of curvature of $S_t=f_t(S^n)\subset{\mathbb R}^{n+1}$. 

As noted in \cite{andrews}, this flow, referred to there as the inverse harmonic mean curvature flow, is expanding and the support function $r$ of the 
surface evolves by the linear strictly parabolic equation
\[
\left(\frac{\partial }{\partial t}-\triangle_{S^n}\right)r=nr.
\]
As a result, the support function for a closed convex surface increases exponentially and the surface expands to infinity. 
Moreover, in  \cite{sm} it is proven that by rescaling the flow about the origin, the surface converges to a round sphere with centre at 0.

In what follows we extract more information about this flow for $n=2$ and prove that

\vspace{0.1in}
\noindent{\bf Main Theorem}:

{\it
Under mean radius of curvature flow, the oriented normal lines to any initial surface with support function $r_0$, converge to those of the round sphere 
with centre at $(x^1,x^2,x^3)$ given by
\begin{align}\label{e:centroid}
x^1+ix^2&={\textstyle{\frac{3}{4\pi}}}\int^{2\pi}_0\int^{\pi}_0r_0(\theta,\phi)\;\sin^2\theta e^{i\phi}d\theta d\phi,\nonumber\\
x^3&={\textstyle{\frac{3}{4\pi}}}\int^{2\pi}_0\int^{\pi}_0r_0(\theta,\phi)\;\sin\theta\cos\theta d\theta d\phi
\end{align},
}
where $(\theta,\phi)$ are standard spherical coordinates. 
\vspace{0.1in}

Thus, these integrals are invariant under the flow.
We prove this by computing the flow of the oriented normal lines to the surface as it evolves by MRCF. These normals form a 
surface in the space of all oriented lines.

In particular, recall that, given any smooth oriented convex surface $S_t$ in 
${\mathbb R}^3$, the set of oriented normal lines to $S_t$ forms a surface $\Sigma_t$ in the space ${\mathbb L}({\mathbb R}^3)$ of all oriented lines 
of ${\mathbb R}^3$. 
This surface is of necessity Lagrangian with respect to the canonical symplectic structure on ${\mathbb L}({\mathbb R}^3)$ and, since $S_t$ is convex, 
$\Sigma_t$ is a section of the bundle $\pi:{\mathbb L}({\mathbb R}^3)=T{\mathbb S}^2\rightarrow {\mathbb S}^2$ \cite{gak2} \cite{gak4}. 

The evolution of the Lagrangian section induced by MRCF, in contrast to the flow in ${\mathbb R}^3$, 
converges without rescaling, and we prove that: 

\vspace{0.1in}
\noindent{\bf Main Theorem} (reformulated):
{\it

Under mean radius of curvature flow, any initial Lagrangian section $F:{\mathbb{C}}\rightarrow {\mathbb{C}}$ converges smoothly to a quadratic 
holomorphic Lagrangian section 
\[
F={\textstyle{\frac{1}{2}}}[x^1+ix^2-2x^3\xi-(x^1-ix^2)\xi^2]
\]
with co-efficients given by the formulae (\ref{e:centroid}).}
\vspace{0.1in}

A holomorphic Lagrangian section corresponds to those oriented lines that pass through a fixed point in ${\mathbb R}^3$, and so, while the 
convex surfaces in ${\mathbb R}^3$ run out to infinity under the flow, their normal lines converge to the normals of a round sphere (without rescaling
the flow). Note that a surface can be both Lagrangian and holomorphic in this setting because the associated Kaehler metric is of neutral signature
\cite{gak4}. 

The proof involves showing that under MRCF, the Lagrangian sections flow by a linear strictly parabolic equation system - see Proposition \ref{p:fflow}. 
Then, utilizing  spherical harmonics to solve the equation in terms of the initial spectral decomposition we study the asymptotic behaviour. 
The fundamental result for the parabolic equation that we use is contained in Proposition \ref{p:lpe}.

In the next section we describe the geometric relationship between ${\mathbb R}^3$ and ${\mathbb L}({\mathbb R}^3)$. Comparison of MRCF in the two
spaces is done in  section 2, while in the final section we prove the Main Theorem.

\vspace{0.1in}


\section{The Space of Oriented Lines}

The space ${\mathbb L}({\mathbb R}^3)$ of oriented lines of Euclidean ${\mathbb R}^3$ can be identified with $T{\mathbb S}^2$, the total space of the 
tangent bundle to the 2-sphere. $T{\mathbb S}^2$ carries a neutral K\"ahler structure ${\mathbb G},{\mathbb J},\Omega)$ which is invariant under the 
Euclidean group acting on oriented lines. In 
what follows, the terms holomorphic and Lagrangian refer to the complex structure ${\mathbb J}$ and symplectic structure $\Omega$, respectively. The metric
${\mathbb G}$ is of neutral signature, hence planes can be both holomorphic and Lagrangian. Further details on the neutral K\"ahler
structure  can be found in \cite{gak2} \cite{gak4}.

For local computations, let $\xi$ be the standard complex coordinate on ${\mathbb S}^2$ coming from stereographic projection from the south pole so that
$\xi=\tan({\textstyle{\frac{\theta}{2}}})e^{i\phi}$ for the spherical polar coordinates $0\leq\theta\leq\pi$ and $0\leq\phi<2\pi$. Extend this
to complex coordinates $(\xi,\eta)$ on an open set of $T{\mathbb S}^2$ by identifying $X\in T_\xi {\mathbb S}^2$ with $(\xi,\eta)\in{\mathbb C}^2$ when
\[
X=\eta\frac{\partial}{\partial\xi}+\bar{\eta}\frac{\partial}{\partial\bar{\xi}}.
\]

Consider the set of oriented normal lines to a surface $S$. These form a Lagrangian surface $\Sigma\subset T{\mathbb S}^2$. As there are no flat points, the 
Gauss map of  $S$ is invertible and hence, $\Sigma$ is a Lagrangian section of the canonical bundle $\pi: T{\mathbb S}^2\rightarrow {\mathbb S}^2$. 
In terms of local 
coordinates the surface $\Sigma$ is given by $\xi\mapsto(\xi,\eta=F(\xi,\bar{\xi}))$ for some complex valued function $F$. 

The link between these holomorphic coordinates and flat coordinates $(x^1,x^2,x^3)$ in ${\mathbb R}^3$ is provided by the map 
$\Phi:T{\mathbb S}^2\times {\mathbb R}\rightarrow {\mathbb R}^3$:
\begin{equation}\label{e:coord1}
x^1+ix^2=\frac{2(\eta-\bar{\eta}\xi^2)+2\xi(1+\xi\bar{\xi})r}{(1+\xi\bar{\xi})^2}
\qquad\qquad
x^3=\frac{- 2(\eta\bar{\xi}+\bar{\eta}\xi)+(1-\xi^2\bar{\xi}^2)r}{(1+\xi\bar{\xi})^2},
\end{equation}
which sends an oriented line $(\xi,\eta)$ and a real number $r$ to the point on the line in ${\mathbb R}^3$ that is an oriented distance $r$ from the 
closest point on the line to the origin.

These equations can be recast as
\[
\eta={\textstyle{\frac{1}{2}}}(x^1+ix^2-2x^3\xi-(x^1-ix^2)\xi^2) \qquad
          r=\frac{(x^1+ix^2)\bar{\xi}+(x^1-ix^2)\xi+x^3(1-\xi\bar{\xi})}{1+\xi\bar{\xi}}.
\]
The perpendicular distance $\chi$ of an oriented line $(\xi,\eta)$ to the origin is found to be
\begin{equation}\label{e:perpdist}
\chi^2=\frac{4\eta\bar{\eta}}{(1+\xi\bar{\xi})^2}.
\end{equation}

\vspace{0.1in}
\begin{Def}
The {\it support function} of  a convex surface is the map $r:S\rightarrow{\mathbb R}$ which takes a point $p$ to the signed distance between 
$p$ and the point on the oriented normal line to $S$ at $p$ which lies closest to the origin. Alternatively it is the signed perpendicular distance 
between the oriented tangent plane to $S$ at $p$ and the origin.
\end{Def}
\vspace{0.1in}

The relationship between the support function $r$ of $S$ and the Lagrangian section $F$ is
\begin{equation}\label{e:supp}
F={\textstyle{\frac{1}{2}}}(1+\xi\bar{\xi})^2\bar{\partial}r.
\end{equation}
Define the complex slopes of $F$ by
\begin{equation}\label{e:slopes}
\bar{\partial}F=-\bar{\sigma} \qquad \qquad (1+\xi\bar{\xi})^2\partial\left(\frac{F}{(1+\xi\bar{\xi})^2}\right)=\rho+i\lambda.
\end{equation}
A section is Lagrangian iff $\lambda=0$ and this implies the existence of the real function $r$ satisfying equation (\ref{e:supp}). In addition,
the radii of curvature of the surface $S$ are determined by
\[
|\sigma|^2={\textstyle{\frac{1}{4}}}(r_1-r_2)^2
\qquad\qquad
(r+\rho)^2={\textstyle{\frac{1}{4}}}(r_1+r_2)^2.
\]
Finally, translations in ${\mathbb R}^3$ act on our functions as follows. Suppose we consider the translation that takes the origin to 
$(x^1+ix^2,x^3)=(\alpha,b)$. Then we have 
\[
\eta\mapsto\eta+{\textstyle{\frac{1}{2}}}(\alpha-2b\xi-\bar{\alpha}\xi^2) \qquad
          r\mapsto r+\frac{\alpha\bar{\xi}+\bar{\alpha}\xi+b(1-\xi\bar{\xi})}{1+\xi\bar{\xi}},
\]
while $\sigma$ and $r+\rho$ are invariant under translations.

\vspace{0.2in}

\section{Mean Radius of Curvature Flow}

Let us now consider the flow (\ref{flow}) for a strictly convex surface $S_t$ in ${\mathbb R}^3$. Using coordinates $(x^1+ix^2,x^3)$ on ${\mathbb R}^3$ 
and Gauss coordinates 
$\xi$ on $S_t$, let $r_t:S^2\rightarrow{\mathbb R}$ be the support function of $S_t$. Then, differentiating
equations (\ref{e:coord1}) in time
\[
\frac{\partial}{\partial t}(x^1+ix^2)=\frac{2}{(1+\xi\bar{\xi})^2}\frac{\partial}{\partial t}\eta
-\frac{2\xi^2}{(1+\xi\bar{\xi})^2}\frac{\partial}{\partial t}\bar{\eta}
+\frac{2\xi}{1+\xi\bar{\xi}}\frac{\partial}{\partial t}r,
\]
\[
\frac{\partial}{\partial t}x^3=-\frac{2\bar{\xi}}{(1+\xi\bar{\xi})^2}\frac{\partial}{\partial t}\eta
-\frac{2\xi}{(1+\xi\bar{\xi})^2}\frac{\partial}{\partial t}\bar{\eta}
+\frac{1-\xi\bar{\xi}}{1+\xi\bar{\xi})}\frac{\partial}{\partial t}r,
\]
and projecting we obtain
\[
\frac{\partial f}{\partial t}^\perp=\frac{\partial r}{\partial t}N=(r_1+r_2)N.
\]
Finally, from the relationship between section and support (\ref{e:supp}) we have
\begin{align}
r_1+r_2&=2(r+\rho)=2r+2(1+\xi\bar{\xi})^2\partial\left(\frac{F}{(1+\xi\bar{\xi})^2}\right)\nonumber\\
&=2r+(1+\xi\bar{\xi})^2\partial\bar{\partial}r =2r+\triangle_{{\mathbb S}^2}r.\nonumber
\end{align}
We have therefore proven the first part of
\vspace{0.1in}
\begin{Prop}
Under MRCF the support function evolves by
\begin{equation}\label{e:suppflow}
\left(\frac{\partial}{\partial t}-\triangle_{{\mathbb S}^2}\right)r=2r,
\end{equation}
while the perpendicular distance function of the normal lines evolves by 
\[
\left(\frac{\partial}{\partial t}-\triangle_{{\mathbb S}^2}\right)\chi^2=2\chi^2-4(\rho^2+|\sigma|^2).
\]
\end{Prop}
\begin{pf}
The first we have proven, the second follows from a similar calculation.
\end{pf}
\vspace{0.1in}

In the space of oriented lines, the set of oriented normal lines to $S_t$ form a Lagrangian section of $T{\mathbb S}^2\rightarrow {\mathbb S}^2$ given locally by a complex function
$F:{\mathbb C}\rightarrow{\mathbb C}$. We lift the flow to the space of oriented lines:

\vspace{0.1in}
\begin{Prop}\label{p:fflow}
Under MRCF the Lagrangian section $F$ evolves in $T{\mathbb S}^2$ by the linear parabolic system
\[
\left(\frac{\partial}{\partial t}-\triangle_{{\mathbb S}^2}\right)F=-\frac{2\bar{\xi}}{1+\xi\bar{\xi}}\bar{\partial}F.
\]
\end{Prop}
\begin{pf}
Differentiate the relationship (\ref{e:supp}) in time to get
\[
\frac{\partial}{\partial t}F={\textstyle{\frac{1}{2}}}(1+\xi\bar{\xi})^2\bar{\partial}\frac{\partial}{\partial t}r
=\triangle_{{\mathbb S}^2}F-\frac{2\bar{\xi}}{1+\xi\bar{\xi}}\bar{\partial}F.
\]
\end{pf}
\vspace{0.1in}

Finally, computing the flow of the derived quantities:
\vspace{0.1in}
\begin{Prop}\label{p:slopesflow}
Under MRCF the slopes evolve by
\[
\left(\frac{\partial}{\partial t}-\triangle_{{\mathbb S}^2}\right)\rho=2\rho
\qquad\qquad
\left(\frac{\partial}{\partial t}-\triangle_{{\mathbb S}^2}\right)\lambda=-2\lambda,
\]
\[
\left(\frac{\partial}{\partial t}-\triangle_{{\mathbb S}^2}\right)\sigma=-2(1+2\xi\bar{\xi})\sigma+
  2(1+\xi\bar{\xi})(\bar{\xi}\bar{\partial}\sigma-\xi\partial\sigma).
\]
\end{Prop}
\begin{pf}
Differentiate the defining relationships (\ref{e:slopes}) in time and use the previous Proposition.
\end{pf}
\vspace{0.1in}

Note that the flow equation for $\lambda$ is such that, if $\lambda=0$ initially, it remains so for all time. Since 
$\Omega|_\Sigma=\lambda d_{\mathbb S}^2A$, we say that the flow in $T{\mathbb S}^2$ is Lagrangian, since it preserves the Lagrangian condition. 

In fact, even if the intial 
surface is not Lagrangian, we prove in Proposition \ref{p:nonlag} that under the flow it becomes Lagrangian.

\vspace{0.2in}

\section{Proof of the Main Theorem}

Consider the flow
\begin{equation}\label{e:spfs2}
\left(\frac{\partial}{\partial t}-\triangle_{{\mathbb S}^2}\right)f=2f,
\end{equation}
for $f:S\times[0,\infty)\rightarrow{\mathbb R}$ with $f(\cdot,0)=f_0(\cdot)$. 

\vspace{0.1in}
\begin{Def}
Define the {\it spherical area} $A_{{\mathbb S}^2}(f)$ of $f$ by:
\[
A_{{\mathbb S}^2}(f)=\int\!\!\!\int_{{\mathbb S}^2} f\;dA.
\]
\end{Def}

\vspace{0.1in}
\begin{Prop}\label{p:lpe}
The above flow converges if and only if the $A_{{\mathbb S}^2}(f_0)=0$. 

For  $A_{{\mathbb S}^2}(f_0)=0$, it converges smoothly to an eigenfunction for the spherical Laplacian 
with eigenvalue 2.

For $A_{{\mathbb S}^2}(f_0)\neq 0$, there exists a constant $C$ depending only $f_0$ such that
\[
|f|\geq Ce^t.
\]
\end{Prop}

\begin{pf}
The flow (\ref{e:spfs2}) is linear and strictly parabolic, and therefore by standard theory \cite{Lieberman}, given any initial function,  there exists a smooth 
solution for all time. Let $f_t$ be the solution of the flow for some initial $f_0$. 

Integrating the flow equation over the 2-sphere 
\[
\frac{\partial}{\partial t}A_{{\mathbb S}^2}(f)=A_{{\mathbb S}^2}(f).
\]
Thus if $A_{{\mathbb S}^2}(f_0)=0$, then $A_{{\mathbb S}^2}(f)=0$ for all time, while $A_{{\mathbb S}^2}(f_0)\neq0$ implies exponential growth in time for the spherical area.

For fixed time $t$, decompose $f_t:S\rightarrow{\mathbb R}$ in terms of spherical harmonics $Y_l^m:S\rightarrow{\mathbb R}$ \cite{muller}:
\[
f_t=\sum_{l=0}^\infty\sum_{m=-l}^{m=l}B_{lm}Y_l^m,
\]
where $B_{lm}$ are complex and satisfy $\overline{B_{lm}}=(-1)^{l}B_{l-m}$ for $m\neq0$ and $\overline{B_{l0}}=B_{l0}$.

Since the flow is linear, we obtain a flow on the projection of $f$ onto the spectrum of the Laplacian:
\[
\frac{\partial B_{lm}}{\partial t}=[2-l(l+1)]B_{lm},
\]
which integrate to yield
\[
f_t=\sum_{l=0}^\infty\sum_{m=-l}^{m=l}\stackrel{\circ}{B}_{lm}e^{[2-l(l+1)]t}Y_l^m.
\]
For convenience we have denoted ${B}_{lm}$ at $t=0$ by $\stackrel{\circ}{B}_{lm}$.

Splitting off the first few terms
\begin{align}
f_t&=\stackrel{\circ}{B}_{00}e^{2t}Y_0^0+\stackrel{\circ}{B}_{1-1}Y_1^{-1}+\stackrel{\circ}{B}_{10}Y_1^0+\stackrel{\circ}{B}_{11}Y_1^1+
\sum_{l=2}^\infty\sum_{m=-l}^{m=l}\stackrel{\circ}{B}_{lm}e^{(2-l(l+1))t}Y_l^m\nonumber\\
&={\textstyle{\frac{1}{2}}}{\textstyle{\sqrt{\frac{1}{\pi}}}}\stackrel{\circ}{B}_{00}e^{2t}
+{\textstyle{\sqrt{\frac{3}{2\pi}}}}\stackrel{\circ}{B}_{1-1}\frac{\bar{\xi}}{1+\xi\bar{\xi}}
+{\textstyle{\frac{1}{2}}}{\textstyle{\sqrt{\frac{1}{\pi}}}}\stackrel{\circ}{B}_{10}\frac{1-\xi\bar{\xi}}{1+\xi\bar{\xi}}
-{\textstyle{\sqrt{\frac{3}{2\pi}}}}\stackrel{\circ}{B}_{11}\frac{\xi}{1+\xi\bar{\xi}}\nonumber\\
&\qquad\qquad\qquad\qquad+\sum_{l=2}^\infty\sum_{m=-l}^{m=l}\stackrel{\circ}{B}_{lm}e^{(2-l(l+1))t}Y_l^m\nonumber\\
&={\textstyle{\frac{1}{2}}}{\textstyle{\sqrt{\frac{1}{\pi}}}}A_{{\mathbb S}^2}(f_0)e^{2t}
-\frac{\alpha\bar{\xi}+\bar{\alpha}\xi+b(1-\xi\bar{\xi})}{1+\xi\bar{\xi}}+\sum_{l=2}^\infty\sum_{m=-l}^{m=l}\stackrel{\circ}{B}_{lm}e^{(2-l(l+1))t}Y_l^m,\nonumber
\end{align}
where we note that, by the orthogonality properties of the spherical harmonics
\[
A_{{\mathbb S}^2}(f_0)=\int\!\!\!\int_{{\mathbb S}^2} f_0\;dA=\stackrel{\circ}{B}_{00},
\]
and we have introduced $\alpha\in{\mathbb C}$,$b\in{\mathbb R}$:
\[
\alpha=-{\sqrt{\textstyle{\frac{3}{2\pi}}}}\stackrel{\circ}{B}_{11}
\qquad\qquad
b={\textstyle{\frac{1}{2}}}{\sqrt{\textstyle{\frac{1}{\pi}}}}\stackrel{\circ}{B}_{10}.
\]

If $A_{{\mathbb S}^2}(f_0)\neq0$, then $|f_t|$ exponentially blows up as $t\rightarrow \infty$, while for $A_{{\mathbb S}^2}(f_0)=0$, $f_t$ converges to an eigenfunction 
of the spherical Laplacian with eigenvalue equal to 2:
\[
f_t\rightarrow -\frac{\alpha\bar{\xi}+\bar{\alpha}\xi+b(1-\xi\bar{\xi})}{1+\xi\bar{\xi}}
\qquad{\mbox{as}}\qquad t\rightarrow \infty,
\]
as claimed.
\end{pf}
\vspace{0.1in}

\noindent{\bf Proof of Main Theorem.}

We have seen that under MRCF the support flows by equation (\ref{e:suppflow}). By Proposition \ref{p:lpe} to determine the behaviour we need
to compute the area $A_{{\mathbb S}^2}(r)$. Note that this integral is invariant under translation:
\[
A_{{\mathbb S}^2}(r)=A_{{\mathbb S}^2}\left(r-\frac{\alpha\bar{\xi}+\bar{\alpha}\xi+b(1-\xi\bar{\xi})}{1+\xi\bar{\xi}}\right).
\] 
If we move the origin to the interior of $S$ we can make $r>0$, and conclude that $A_{{\mathbb S}^2}(r)>0$. Thus, by Proposition \ref{p:lpe},
the support function blows up exponentially. More particularly, the spectral decomposition is
\[
r=r_{00}e^{2t}+\sum_{l=1}^\infty\sum_{m=-l}^{m=l}r_{lm}e^{[2-l(l+1)]t}Y_l^m,
\] 
where $r_{00}>0$. Clearly, $r_{00}$ is the radius of the limit of the rescaled flow for $\tilde{r}=re^{-2t}$ \cite{sm} - it is the projection of 
the support function into the 0-eigenspace of the spherical Laplacian:
\[
r_{00}={\textstyle{\frac{1}{4\pi}}}\int^{2\pi}_0\int^{\pi}_0r_0\sin\theta d\theta d\phi
\]
Here we have introduced standard spherical polar coordinates $\xi=\tan({\textstyle{\frac{\theta}{2}}})e^{i\phi}$.

Now consider the flow in $T{\mathbb S}^2$. By the evolution equation for $\rho$ in Proposition \ref{p:slopesflow} we require
$A_{{\mathbb S}^2}(\rho)$. In this case, however we have
\[
A_{{\mathbb S}^2}(\rho)={\textstyle{\frac{1}{2}}}\int\!\!\int_{{\mathbb S}^2}\Delta_{{\mathbb S}^2}r\;dA=0,
\]
and by Proposition \ref{p:lpe}
\[
\rho=\sum_{l=1}^\infty\sum_{m=-l}^{m=l}\rho_{lm}e^{(2-l(l+1))t}Y_l^m,
\]
where the constants $\rho_{lm}$ are determined by the initial surface. 
Let
\[
\alpha=-{\sqrt{\textstyle{\frac{3}{2\pi}}}}\rho_{11}
\qquad\qquad
b={\textstyle{\frac{1}{2}}}{\sqrt{\textstyle{\frac{1}{\pi}}}}\rho_{10}.
\]
We claim that the oriented normal lines to the flowing surface converge to the set of oriented lines passing through $(x^1+ix^2,x^3)=(\alpha,b)$.
Using the power series expansion for $r$ corresponding to the one for $f_t$ in the proof of Proposition 4, 
we have 
\[F = {\textstyle{\frac{1}{2}}}(1 + \xi \bar\xi)^2\bar\partial r =  
\frac{1}{2}(\alpha - 2b\xi - \bar\alpha \xi^2) + R(t,\xi,\bar\xi),
 \]
where $ R(t,\xi,\bar\xi) \to 0$ uniformly in any $C^k({\mathbb S}^2)$ as $t \to \infty$. This completes the proof of the Main 
Theorem since the first term represents a global Lagrangian holomorphic section and we have
\begin{align}
x^1+ix^2&=\alpha=-{\sqrt{\textstyle{\frac{3}{2\pi}}}}\rho_{11}=-{\sqrt{\textstyle{\frac{3}{2\pi}}}}\int\!\!\int_{{\mathbb S}^2}Y_1^1\rho \;dA\nonumber\\
               &=-{\textstyle{\frac{1}{2}}}{\sqrt{\textstyle{\frac{3}{2\pi}}}}\int\!\!\int_{{\mathbb S}^2}Y_1^1\triangle_{{\mathbb S}^2}r \;dA\nonumber\\
               &={\sqrt{\textstyle{\frac{3}{2\pi}}}}\int\!\!\int_{{\mathbb S}^2} Y_1^1 r\;dA\nonumber\\
               &={\textstyle{\frac{3}{4\pi}}}\int^{2\pi}_0\int^{\pi}_0r\sin^2\theta e^{i\phi}\;d\theta d\phi\nonumber,
\end{align}
and similarly for $x^3$.
\qed
\vspace{0.2in}

In fact, one could drop the Lagrangian condition on the initial global section:

\vspace{0.2in}
\begin{Prop}\label{p:nonlag}
The flow converges to a Lagrangian holomorphic section even if the initial section is not Lagrangian.
\end{Prop}
\begin{pf}
Recall that the section is Lagrangian iff $\lambda=0$. The flow equation for $\lambda$ given in Proposition \ref{p:slopesflow} can be written
\[
\left(\frac{\partial }{\partial t}-\triangle_{{\mathbb S}^2}\right)\lambda e^{2t}=0,
\]
which by the maximum principle implies that there exists a constant $C$ depending only on $\lambda_0$ such that
\[
|\lambda|\leq Ce^{-2t},
\]
so that $\lambda\rightarrow 0$ as $t\rightarrow \infty$. 

Now the proof follows that of the Lagrangian case.
\end{pf}

\end{document}